\documentclass[preprints,article,accept,moreauthors,pdftex]{Definitions/mdpi} 

\firstpage{1} 
\pubvolume{xx}
\issuenum{1}
\articlenumber{5}
\pubyear{2019}
\copyrightyear{2019}
\history{}

\pdfoutput=2

\usepackage{amsfonts}
\usepackage{amssymb}
\usepackage{graphicx}
\usepackage{hyperref,color}
\usepackage{comment}%
\providecommand{\U}[1]{\protect\rule{.1in}{.1in}}

\newcommand{\rem}[1]{}

\numberwithin{equation}{section}

\Title{Determining When An Algebra Is An Evolution Algebra}

\Author{Miguel D. Bustamante$^{1,\dagger}$\orcidA{}, Pauline Mellon $^{2,\dagger}$\orcidB{} and M. Victoria Velasco $^{3,\dagger}$*\orcidC{}}

\AuthorNames{Miguel D. Bustamante, Pauline Mellon and M. Victoria Velasco}

\address{%
$^{1}$ \quad School of Mathematics and Statistics, University College Dublin, Dublin 4, Ireland; miguel.bustamante@ucd.ie\\
$^{2}$ \quad School of Mathematics and Statistics, University College Dublin, Dublin 4, Ireland;
pauline.mellon@ucd.ie\\
$^{3}$ \quad Departamento de An\'{a}lisis Matem\'{a}tico,
Facultad de Ciencias, Universidad de Granada, 18071 Granada, Spain}

\corres{Correspondence: vvelasco@ugr.es; Tel.: +34-958-241-000}

\firstnote{These authors contributed equally to this work.}

\abstract{Evolution algebras are non-associative algebras that describe non-Mendelian
hereditary processes and have connections with many other areas.
In this paper we obtain necessary and sufficient conditions for a given
algebra $A$ to be an evolution algebra. We prove that the problem is
equivalent to the so-called \emph{SDC problem,} that is, the\emph{
simultaneous diagonalisation via congruence} of a given set of matrices\emph{.
}More precisely we show that an $n$-dimensional algebra $A$ is an evolution
algebra if, and only if, a certain set of $n$ symmetric $n\times n$ matrices
$\{M_{1},\ldots,M_{n}\}$ describing the product of $A$ are \emph{SDC}. We
apply this characterisation to show that while certain classical genetic
algebras (representing Mendelian and auto-tetraploid inheritance) are not
themselves evolution algebras, arbitrarily small perturbations of these are
evolution algebras. This is intriguing as evolution algebras model asexual
reproduction unlike the classical ones.}


\begin{document}


\section{Introduction}

This is a postprint of our work published in: \emph{Mathematics} \textbf{2020}, \emph{8}, 1349. \url{https://doi.org/10.3390/math8081349}. It contains a few minor modifications on pages 2 and 4.

Evolution algebras are non-associative algebras with a dynamic nature. They
were introduced in 2008 by Tian \cite{Tian} to enlighten the study of
non-Mendelian genetics. Since then, a large literature has flourished on this
topic (see for instance \cite{Be-Me-Ve,
Ca-Si-Ve,Ce-Ve,Ca-Go,3,5,Elduque,He,13,15,Mellon-Ve, Ro-Ve,19,20,21,Ve})
motivated by the fact that these algebras have connections with group theory,
Markov processes, theory of knots, systems and graph theory. For instance, in
\cite{Be-Me-Ve}, the theory of evolution algebras was related to that of pulse
processes on weighted digraphs and applications were provided by reviewing a report of the National Science Foundation about air pollution
achieved by the Rand Corporation. A pulse process is a structural dynamic
model to analyse complex networks by studying the propagation of changes
through the vertices of a weighted digraph, after introducing an initial pulse
in the system at a particular vertex. It is based on a spectral analysis of
the corresponding weighted digraph to facilitate large scale decision making
processes. Evolution algebras also introduce useful
algebraic techniques into the study of some digraphs because evolution
algebras and weighted digraphs can be canonically identified.

We recall that an \textbf{algebra} is a linear space $A$ provided with a
product, that is, a bilinear map from $A\times A$ to $A$ via the operation
$(a,b)\rightarrow ab$. In the particular case that $(ab)c=a(bc),\,\,\text{for
all}\,\,a,b,c\in A$ we say that $A$ is\textbf{ }associative. Meanwhile, if
$ab=ba,\,\,\text{for all}\,\,a,b\in A,$ then we say that $A$ is commutative.

An \textbf{evolution algebra }is defined as a commutative algebra $A$ for
which there exists a basis $B^{\ast}=\{e_{i}^{\ast}:i\in\Lambda\}$ such that
$e_{i}^{\ast}e_{j}^{\ast}=0$ for every $i,j\in\Lambda$ with $i\neq j$. Such a
basis is called\textbf{ natural}. Evolution algebras are, in general,
non-associative. To date most literature on evolution algebras is on
finite-dimensional ones. However, in \cite{Mellon-Ve} it is shown that every
infinite-dimensional Banach evolution algebra is the direct sum of a
finite-dimensional evolution algebra and a zero-product algebra.

In this paper we discuss necessary and sufficient conditions under which a
given finite-dimensional commutative algebra is an evolution algebra, namely,
we determine when such a finite-dimensional algebra can be provided with a
natural basis. We tackle the problem constructively by assuming an arbitrary
basis $B$ with a multiplication table given by equation (\ref{eq:mult_table})
below and then asking whether or not there is a change of basis from $B$ to a
natural basis $B^{\ast}.$ In Section 2, Theorem \ref{thm-main}, we show that
this problem is equivalent to the simultaneous diagonalisation via congruence
of certain $n\times n$ symmetric matrices $M_{1},\ldots,M_{n},$ called the
multiplication structure matrices obtained from the given multiplication table.

Finding concrete sufficient conditions for a given set of matrices to be
simultaneously diagonalisable via congruence (we will refer to it as the
SDC-problem) is one of the 14 open problems posted in 1990 by Hiriart-Urruty
\cite{Hi-U} (see also \cite{Hi-Ur-Tor, Ji-Li}). It has connections with other
problems such as blind-source separation in signal processing
\cite{belouchrani1997blind, yeredor2000blind, yeredor2002non,Vol}. The
SDC-problem was solved recently for complex symmetric matrices in
\cite{Bu-Me-Ve-SD}.

In Theorem \ref{complex1} we show that if $A$ is a real algebra and $B$ is a
basis of $A$ then $B\,$\ also is a basis of \thinspace$A_{\mathbb{C}},$ the
complexification of $A$ (with the same multiplication structure matrices) and
that $A$ is an evolution algebra if, and only if, \thinspace$A_{\mathbb{C}}$
is an evolution algebra {and has a natural basis consisting of elements of $A$}. This reduction of the real case to the complex one
allows us to apply the results in \cite{Bu-Me-Ve-SD} to both real and complex algebras.

In Theorem \ref{rankn}\ we determine if a given algebra $A$ whose annihilator
is zero is an evolution algebra and in Theorem \ref{dos} we do the same if its
annihilator is not zero. A useful characterisation of the property of being an
evolution algebra is given in the particular case that one of the
multiplication structure matrices is invertible. In this case if $M_{i_{0}}$
is invertible then $A$ is an evolution algebra if, and only if, for each
$k\neq i_{0}$ the matrix $M_{i_{0}}^{-1}M_{k}\,\ $is diagonalisable by
similarity and these matrices pairwise commute.

Applications of these results are provided in the final section of this paper.
They also show that the conditions in the mentioned results are neither
redundant nor superfluous.

We prove that some classical genetic algebras such as the gametic algebra for
simple Mendelian inheritance (Example \ref{ejem1}) or the gametic algebra for
auto-tetraploid inheritance (Example \ref{ejem2}) are not evolution algebras.
Nevertheless, both of these algebras can be deformed by means of a parameter
$\varepsilon>0$ to obtain an algebra $A_{\varepsilon}$ that is an evolution
algebra for every value of the parameter $\varepsilon,$ as shown in Example
\ref{ejem1de} and Example \ref{ejem2de} respectively.

\section{Characterising evolution algebras by means of simultaneous
diagonalisation of matrices by congruence}

An $n$-dimensional algebra $A$ over a field $\mathbb{K}$ ($\mathbb{=R}$ or
$\mathbb{C}$) is determined by means of a basis $B=\{e_{1},\ldots,e_{n}\}$
together with a multiplication table

\begin{equation}
e_{i}e_{j}={\sum\limits_{k=1}^{n}}m_{ijk}\,e_{k}\,,\quad i,j=1,\ldots,n\,,
\label{eq:mult_table}%
\end{equation}
where $m_{ijk}\in\mathbb{K}\,,$ for $i,j,k=1,\ldots,n$. \ In fact, if
$a:=\sum\limits_{i=1}^{n}\alpha_{i}e_{i}$ and $b:=\sum\limits_{j=1}^{n}%
\beta_{j}e_{j}$ then, by bilinearity, the product $ab$ is obtained from the
multiplication table (\ref{eq:mult_table}) as follows
\[
ab=\left(  \sum\limits_{i=1}^{n}\alpha_{i}e_{i}\right)  \left(  \sum
\limits_{j=1}^{n}\beta_{j}e_{j}\right)  =\sum\limits_{k=1}^{n}\left(
\sum\limits_{i,j=1}^{n}\alpha_{i}\beta_{j}m_{ijk}\right)  e_{k},
\]
where $m_{ijk}:=\pi_{k}(e_{i}e_{j})$ and $\pi_{k}:A\rightarrow\mathbb{K}$ is
the projection over the $k$-th coordinate, that is $\pi_{k}(\sum
\limits_{i=1}^{n}\alpha_{i}e_{i})=\alpha_{k}.$

These basis-dependent coefficients $m_{ijk}$ are known as structure
constants\textbf{ }with respect to $B$ (see \cite{Pierce82}). For a basis $B$
of $A$, the structure constants completely determine the algebra $A$, up to
isomorphism. \ 

If we organise the $n^{3}$ \ structure constants in $n$ matrices by defining

\begin{equation}
M_{k}(B):=\left(
\begin{array}
[c]{ccc}%
\pi_{k}(e_{1}e_{1}) &  & \pi_{k}(e_{1}e_{n})\\
\vdots &  & \vdots\\
\pi_{k}(e_{n}e_{1}) &  & \pi_{k}(e_{n}e_{n})
\end{array}
\right)  =\left(
\begin{array}
[c]{ccc}%
m_{11k} &  & m_{1nk}\\
\vdots &  & \vdots\\
m_{n1k} &  & m_{nnk}%
\end{array}
\right)  , \label{matriz}%
\end{equation}
for $k=1,...,n,$ then the product of $A$ is given by%

\begin{equation}
\left(  \sum\limits_{i=1}^{n}\alpha_{i}e_{i}\right)  \left(  \sum
\limits_{j=1}^{n}\beta_{j}e_{j}\right)  =\sum\limits_{k=1}^{n}\left(
\mbox{\boldmath{$\alpha$}}^{T}M_{k}(B)\mbox{\boldmath{$\beta$}}\right)
e_{k}\,, \label{eq:M_k}%
\end{equation}
where $\mbox{\boldmath{$\alpha$}}^{T}=(\alpha_{1},\ldots,\alpha_{n}),$
$\mbox{\boldmath{$\beta$}}^{T}=(\beta_{1},\ldots,\beta_{n})$ and $T$ indicates
the transpose operation. This motivates the following definition.

\begin{Definition}
\label{struct}\emph{If }$A$\emph{ is an algebra, }$B=\{e_{1},\ldots,e_{n}%
\}$\emph{ is a basis of }$A$\emph{ and }$e_{i}e_{j}=\sum\limits_{k=1}%
^{n}m_{ijk}\,e_{k},$\emph{ for }$i,j=1,\ldots,n\,,$\emph{ then the
\textbf{multiplication structure matrices (m-structure matrices }for
short\textbf{)} of }$A$\emph{ with respect to }$B$\emph{ are the }$n\times n$
\emph{matrices }$M_{k}(B)=\left(  \pi_{k}(e_{i}e_{j})\right)  $ \ \emph{given
by (\ref{matriz}) \ for }$k=1,...,n.$\emph{ Note that these matrices are
symmetric if, and only if, }$A\,$\emph{is commutative.}

\emph{If the basis }$B$\emph{ is clear from the context then we will write
}$M_{k}:=M_{k}(B)\,\ $\emph{\ for }$k=1,...,n.$
\end{Definition}

We recall that an $n$-dimensional \textbf{evolution algebra }is a commutative
algebra $A$ for which there exists a basis $B^{\ast}=\{e_{1}^{\ast}%
,...,e_{n}^{\ast}\}$ such that $e_{i}^{\ast}e_{j}^{\ast}=0$ for every
$i,j\in\{1,\cdots,n\}$ with $i\neq j$. Such a basis $B^{\ast}$ is said to be a
\textbf{natural basis} of $A$.

The
next result is a straightforward combination of the concept of evolution algebra with Definition \ref{struct}. 
\begin{Proposition}
\label{anterior}An evolution algebra is an algebra $A$ provided with a basis
$B^{\ast}=\{e_{1}^{\ast},...,e_{n}^{\ast}\}$ such that the corresponding
m-structure matrices $M_{1}(B^{\ast})=(\pi_{1}(e_{i}^{\ast}e_{j}^{\ast})),$
$\cdots,$ $M_{n}(B^{\ast})=(\pi_{n}(e_{i}^{\ast}e_{j}^{\ast}))$ are diagonal.
\end{Proposition}

\begin{proof}
$M_{k}(B^{\ast})$ is diagonal for $k=1,...,n,$ if, and only if, $e_{i}^{\ast
}e_{j}^{\ast}=0,$ for every $i\neq j,$ or equivalently if $B^{\ast}$ is a
natural basis.
\end{proof}
In the next theorem we characterise when a given algebra is an evolution
algebra. To this end we recall the following property.

\begin{Definition}
\label{SDC}\emph{Let }$M_{1},\ldots,M_{m}$\emph{ be symmetric
}$n\times n$\emph{ matrices. Then these matrices are\textbf{ simultaneously
diagonalisable via congruence} (SDC) if, and only if, there exists a
nonsingular }$n\times n$\emph{ matrix }$P$\emph{ and }$m$\emph{ diagonal
}$n\times n$\emph{ matrices }$\{D_{j}\}_{j=1}^{m}$\emph{ such that }%
\[
P^{T}M_{j}P=D_{j},\quad j=1,\ldots,m.
\]

\end{Definition}

It is worth remarking at this point that the general problem of
diagonalisation via congruence considers $m$ symmetric matrices of dimension
$n\times n$, where $m$ need not be equal to $n$. This problem has applications
in statistical signal processing and multivariate statistics
\cite{belouchrani1997blind, yeredor2000blind, yeredor2002non, Vol} and was
solved for complex symmetric matrices in \cite{Bu-Me-Ve-SD}.

\begin{Theorem}
\label{thm-main} Let $A$ be a commutative algebra over $\mathbb{K}$ with basis
$B=\{e_{1},\ldots,e_{n}\}$. $\ $Let $M_{1},\ldots,M_{n}$ be the
m-structure matrices of $A$ with respect to $B.$ Then $A$ is an evolution
algebra if, and only if, the symmetric matrices $M_{1},$ $\ldots,M_{n}$
are simultaneously diagonalisable via congruence.
\end{Theorem}

\begin{proof}
$A$ is an evolution algebra if, and only if, $A$ has a natural basis, say
$B^{\ast}=\{e_{1}^{\ast},...,e_{n}^{\ast}\}$ (that is a basis such that
$e_{i}^{\ast}e_{j}^{\ast}=0$ if $i\neq j$). \ Let $P=(p_{ij})$ be the change
of basis matrix from $B$ to $B^{\ast}$ (that is $e_{i}^{\ast}=\sum_{k=1}%
^{n}p_{ki}{e}_{k}\,$\ for $i=1,\ldots,n$). Then, by (\ref{eq:M_k}),

\begin{equation}
e_{i}^{\ast}e_{j}^{\ast}=\left(  \sum\limits_{k=1}^{n}p_{ki}{e}_{k}\right)
\left(  \sum\limits_{k=1}^{n}p_{kj}{e}_{k}\right)  =\sum\limits_{k=1}%
^{n}\left(  \mbox{\boldmath{$\alpha$}}^{T}M_{k}%
\mbox{\boldmath{$\beta$}}\right)  e_{k}, \label{parti}%
\end{equation}
where $\mbox{\boldmath{$\alpha$}}=P\gamma_{i}$ and
$\mbox{\boldmath{$\beta$}}=P\gamma_{j}$ with $\gamma_{i}%
=(0,...,0,\overset{(i-\text{th})}{1},0,...0)^{T}\in\mathcal{M}_{n\times
1}(\mathbb{K})$. Thus%

\begin{equation}
e_{i}^{\ast}e_{j}^{\ast}=\sum\limits_{k=1}^{n}\left(  \gamma_{i}^{T}P^{T}%
M_{k}P\gamma_{j}\right)  e_{k}=0,\text{ for }i\neq j, \label{parti2}%
\end{equation}
and hence $e_{i}^{\ast}e_{j}^{\ast}=0$ if $i\neq j$ if, and only if, the
matrix $P^{T}M_{k}P$ is diagonal for $k=1,...,n.$
\end{proof}
Since the problem of simultaneous diagonalisation of matrices via congruence
was solved in \cite{Bu-Me-Ve-SD} for complex symmetric matrices, we consider
the following.

The \textbf{complexification }of a real algebra $A$\ is defined as the complex
algebra $\ A_{\mathbb{C}}:=A\oplus iA=\{a+ib:a,b\in A\},$ where for
$a,b,c,d\in A$ and $r,s\in\mathbb{R}$,%

\begin{align*}
(a+ib)+(c+id)  &  =(a+b)+i(b+d),\\
(r+is)(a+ib)  &  =ra-sb+i(rb+sa),\\
(a+ib)(c+id)  &  =(ac-bd)+i(ad+bc).
\end{align*}

Note that every basis $B$ of $A$ is trivially a basis of $A_{\mathbb{C}}$ so
that \ the real dimension of $A$ and the complex dimension of $A_{\mathbb{C}}$ coincide.

\begin{Theorem}
\label{complex1}Let $A$ be a real algebra. Then $A$ is an evolution algebra
if, and only if, $A_{\mathbb{C}}$ is an evolution algebra {and has a natural basis consisting of elements of $A$}. Moreover, if $A$ is
a real evolution algebra then every natural basis of $A$ is a natural basis of
$A_{\mathbb{C}}.$
\end{Theorem}

\begin{proof}
If $A$ is an evolution algebra
and if $B$ is a natural basis of $A$ then obviously $B$ is a natural basis of
$A_{\mathbb{C}}.$ {The converse direction is clear.} 
\end{proof}

\begin{Corollary}
\label{complex2}Let $A$ be a real commutative algebra, \thinspace
$B=\{e_{1},...,e_{n}\}$ a basis and $M_{1},\ldots,M_{n}$ be the
m-structure matrices of $A$ with respect to $B.$ Then $A$ is an evolution
algebra if, and only if, the matrices $M_{1},\ldots,M_{n}$ (regarded as
complex matrices) are simultaneously diagonalisable via congruence {by means of a real matrix}.
\end{Corollary}

{In \cite{Bu-Me-Ve-SD}, example 16, we give two real matrices which are diagonalisable via congruence by means of a complex matrix but not by means of any real matrix.}

\subsection{Reviewing the solution of the SDC problem}

The aim of this subsection is to review the solution of the SDC\ problem, that
is, determining when $m$ matrices of size $n\times n$ are simultaneously
diagonalisable via congruence, which was solved in \cite{Bu-Me-Ve-SD} for
complex matrices. All matrices considered in this section are complex.

From now on, let $\mathcal{M}_{n}$ denote the set of all complex $n\times n$
matrices. Moreover, let $\mathcal{MS}_{n}$ be the set of all symmetric
matrices in $\mathcal{M}_{n}$ and $\mathcal{GL}_{n}$ be the set of nonsingular
matrices in $\mathcal{M}_{n}.$

We recall the following definition of simultaneous diagonalisation of matrices
via similarity (SDS), not to be confused with Definition \ref{SDC} involving
simultaneous diagonalisation via congruence (SDC). Nevertheless, the solution
of the problem of determining when a set of complex matrices is SDC given in
\cite{Bu-Me-Ve-SD} is related to the problem of determining whether a
certain set of related matrices is SDS, as we will show below.

\begin{Definition}
\emph{Let }$N_{1},...,N_{m}\in\mathcal{M}_{n}.$ \emph{These matrices are said
to be \textbf{simultaneously diagonalisable by similarity }(SDS) if, and only
if, there exists }$P\in\mathcal{GL}_{n}$\emph{ }$\ $\emph{such that }%
$P^{-1}N_{k}P$\emph{ is diagonal for every }$k=1,...,m.$
\end{Definition}

The following result is well known \cite[Theorem 1.3.12 and Theorem
1.3.21]{Horn}.

\begin{Proposition}
\label{Horn}Let $N_{1},...,N_{m}\in\mathcal{M}_{n}.$ These matrices are
simultaneously diagonalisable by similarity (SDS) if, and only if, they are
each diagonalisable by similarity and they pairwise commute.
\end{Proposition}

\begin{Remark}
\label{mista}\emph{Concerning the statement of the above theorem in
\cite{Horn} we point out that the fact that the symmetric matrices }$N_{{1}%
},...,N_{m}$\emph{ commute guarantees that }$N_{1},...,N_{m}$ \emph{are
simultaneously diagonalisable by similarity only when each of }$N_{1},...,N_{m}$
\emph{are diagonalisable matrices (and obviously not otherwise).}
\end{Remark}

In \cite{Bu-Me-Ve-SD}, to solve the SDC\ problem, Theorem \ref{rank-n} and
Theorem \ref{rankmenor} below were proved. To state them, we recall the next definition.

\begin{Definition}
\emph{Given }$M_{1},...,M_{m}\in{\mathcal{M}}_{n},$\emph{ define the
associated \textbf{linear pencil} to be the map }$M:\mathbb{C}^{m}$\emph{
}$\rightarrow{\mathcal{M}}_{n}$\emph{ given by }$M(\lambda):=%
{\displaystyle\sum\limits_{j=1}^{m}}
\lambda_{j}M_{j},$\emph{ for every }$\lambda=(\lambda_{1},...,\lambda_{m})$
\emph{in} $\mathbb{C}^{m}.$\emph{ Since, for }$\lambda\neq0$\emph{, }%
\[
{\mathrm{rank} } \,M(\lambda)={\mathrm{rank} } \,M\left(  \frac{\lambda
}{\left\Vert \lambda\right\Vert }\right)  \text{,}%
\]
$\,\ $\ \emph{it follows that }%
\[
\sup\{\mathrm{rank}M(\lambda):\lambda\in\mathbb{C}^{m}\}=\sup\{\mathrm{rank}%
M(\lambda):\lambda\in\mathbb{C}^{m}\text{ \emph{with }}\left\Vert
\lambda\right\Vert =1\}\in\{0,1,...,n\}.
\]
\emph{Consequently, this supremum must be achieved so that there exists
}$\lambda_{0}\in\mathbb{C}^{m}$\emph{ with }$\left\Vert \lambda_{0}\right\Vert
=1$\emph{ such that }%
\[
r_{0}:=\mathrm{rank}\,M(\lambda_{0})=\max\{\mathrm{rank}\,M(\lambda
):\lambda\in\mathbb{C}^{m}\},
\]
\emph{and we say that }$r_{0}$\emph{ is the \textbf{maximum pencil rank }of
}$M_{1},$\emph{ }$...,$\emph{ }$M_{m}.$
\end{Definition}

The next theorem corresponds to Theorem 7 in \cite{Bu-Me-Ve-SD} and deals with
the case when the maximum pencil rank of the matrices is $n.$

\begin{Theorem}
\label{rank-n}Let $M_{1},...,M_{m}\in\mathcal{MS}_{n}$ have maximum pencil
rank $n.$ Let $\lambda_{0}\in\mathbb{C}^{m}$ be such that $r_{0}:=$
\textrm{rank} $M(\lambda_{0})=n.$ Then $M_{1},...,M_{m}$ are SDC if, and only
if, $M(\lambda_{0})^{-1}M_{1},\ ...,M(\lambda_{0})^{-1}M_{m}$ are SDS.
\end{Theorem}

Proposition \ref{Horn} \ gives the following result.

\begin{Corollary}
\label{rankncoro}Let $M_{1},...,M_{m}\in\mathcal{MS}_{n},$ and $\lambda_{0}%
\in\mathbb{C}^{m}$ be such that
\[
r_{0}:=\mathrm{rank}M(\lambda_{0})=n.
\]
Then $M_{1},...,M_{m}$ are SDC if, and only if,$\ M(\lambda_{0})^{-1}%
M_{1},...,M(\lambda_{0})^{-1}M_{m}$ are all diagonalisable by similarity and pairwise
commute$.$
\end{Corollary}

Given $1\leq r<n$, and matrices $M_{r}\in\mathcal{M}_{r}$ and $N_{n-r}%
\in\mathcal{M}_{n-r}$, denote by $M_{r}\oplus N_{n-r}$ the $n\times n$ matrix
given by
\[
\left(
\begin{array}
[c]{cc}%
M_{r} & 0_{r\times(n-r)}\\
0_{(n-r)\times r} & N_{n-r}%
\end{array}
\right)  .
\]

When the pencil rank of $M_{1},...,M_{m}\in\mathcal{MS}_{n}$ is strictly less
than $n\,\ $then the SDC problem can be reduced to a similar one in a reduced
dimension as the following result (Theorem 9 in \cite{Bu-Me-Ve-SD}) shows.

\begin{Theorem}
\label{rankmenor}Let $M_{1},...,M_{m}\in\mathcal{MS}_{n}$ have maximum pencil
rank $r.$ Then the following assertions are equivalent:

$\mathrm{(i)}$ $M_{1},...,M_{m}$ are SDC;

$\mathrm{(ii)}$ $\dim(\cap_{j=1}^{m}\ker M_{j})=n-r$ and there exists
$P\in\mathcal{GL}_{n}$ satisfying $P^{T}M_{j}P=\widetilde{D}_{j}\oplus0_{n-r}$
where $\widetilde{D}_{j}\in\mathcal{MS}_{r}$ is \ diagonal for $1\leq j\leq
m.$

Moreover, if either of the above conditions is satisfied, then the pencil
$\widetilde{D}$ \ associated with the $r\times r$ matrices $\widetilde{D}_{1}%
,...,\widetilde{D}_{m}$ is non-singular. Indeed, if $\lambda_{0}\in
\mathbb{C}^{m}$ with $\left\Vert \lambda_{0}\right\Vert =1$ is such that
$r=\mathrm{rank}\,M(\lambda_{0})$ then $\widetilde{D}(\lambda_{0}%
)\in\mathcal{GL}_{r}.$
\end{Theorem}

\subsection{Checking when an algebra is an evolution algebra}

We apply the above results to the m-structure matrices $M_{1},\ldots,M_{n}$ of
an algebra $A$ with respect to a basis $B=\{e_{1},...,e_{n}\}$ as in
(\ref{matriz}). For a real algebra $A$ we consider the complexification
$A_{\mathbb{C}}$ provided with the same basis $B.$

We recall that the annihilator of an algebra $A$ is the set
\[
\mathrm{Ann}(A)=\{b\in A:ab=ba=0,\text{ for every }a\in A\}.
\]
This set is an ideal of $A.$

\begin{Lemma}
\label{one}Let $A$ be a commutative algebra and $B=\{e_{1},...,e_{n}\}$ be a
basis of $A.$ Let $M_{1},\ldots,M_{n}$ be the m-structure matrices of $A$
with respect to $B$. Then
\[
\mathrm{Ann}(A)=\{%
{\displaystyle\sum_{i=1}^{n}}
\beta_{i}e_{i}:(\beta_{1},...,\beta_{n})^{T}\in\cap_{j=1}^{n}\ker M_{j}\}.
\]

\end{Lemma}

\begin{proof}
Since $\left(  \sum\limits_{i=1}^{n}\alpha_{i}e_{i}\right)  \left(
\sum\limits_{j=1}^{n}\beta_{j}e_{j}\right)  =\sum\limits_{k=1}^{n}\left(
\mbox{\boldmath{$\alpha$}}^{T}M_{k}\mbox{\boldmath{$\beta$}}\right)  e_{k}%
,\,$\ as shown in (\ref{eq:M_k}) we have that if $(\beta_{1},...,\beta
_{n})^{T}\in\cap_{j=1}^{n}\ker M_{j}$ then$\,b:=$ $\sum\limits_{j=1}^{n}%
\beta_{j}e_{j}\in\mathrm{Ann}(A)$ as $ab=ba=0$ for every $a\in A$ (because
$M_{k}\mbox{\boldmath{$\beta$}}=\mathbf{0}$).

Conversely, if $b:=$ $\sum\limits_{j=1}^{n}\beta_{j}e_{j}\in\mathrm{Ann}(A)$
then $e_{i}b=0$ for every $i=1,...,n$.\thinspace\ It follows that,%
\[
(0,...,0,\overset{(i-\text{th})}{1},0,...,0)M_{k}(\beta_{1},...,\beta_{n}%
)^{T}=0,\text{ }%
\]
for $i,k\in\{1,...,n\}.$ Fixing $k$ and running $i$ we deduce that, for each
$k=1,...,n,$%
\[
(\beta_{1},...,\beta_{n})^{T}\in\ker M_{k},
\]
Consequently, $(\beta_{1},...,\beta_{n})^{T}\in\cap_{j=1}^{n}\ker M_{j},$ as desired.
\end{proof}

\begin{Theorem}
\label{rankn}Let $A$ be a complex commutative algebra with $\mathrm{Ann}%
(A)=\{0\}.$ Let $B=\{e_{1},\ldots,e_{n}\}$ be a basis of $A$ and let
$M_{1},\ldots,M_{n}$ be the m-structure matrices of $A$ with respect to $B$.

$\mathrm{(i)}$ If $M_{1},...,M_{n}$ have maximum pencil rank $n$ and
$\lambda_{0}\in\mathbb{C}^{n}$ with $\left\Vert \lambda_{0}\right\Vert =1$ is
such that $\mathrm{rank}$ $M(\lambda_{0})=n$ then $A$ is an evolution algebra
if, and only if, each of the matrices $M(\lambda_{0})^{-1}M_{1},...,M(\lambda
_{0})^{-1}M_{n}$ is diagonalisable by similarity and they pairwise commute.

$\mathrm{(ii)}$\ If $M_{1},\ldots,M_{n}$ have maximum pencil rank $r<n$ then
$A$ is not an evolution algebra.
\end{Theorem}

\begin{proof}
(i) If $\lambda_{0}\in\mathbb{C}^{n}$ with $\left\Vert \lambda_{0}\right\Vert
=1$ is such that $\mathrm{rank}$ $M(\lambda_{0})=n$ then, by Corollary
\ref{rankncoro}, we conclude that $A$ is an evolution algebra if, and only if,
the matrices $M(\lambda_{0})^{-1}M_{1},...,M(\lambda_{0})^{-1}M_{n}$ are
diagonalisable and they pairwise commute. (ii) Otherwise the maximum pencil
rank of $\{M_{1},\ldots,M_{n}\}$ is $r<n$ and, by the above lemma,
$\dim\mathrm{Ann}(A)=\cap_{j=1}^{n}\ker M_{j}=0\neq n-r$. Consequently, by
Theorem \ref{rankmenor}, we conclude that $A$ is not an evolution algebra.
\end{proof}

\begin{Corollary}
\label{diago}Let $A$ be a complex commutative algebra and let $B=\{e_{1}%
,\ldots,e_{n}\}$ be a basis of $A$. Let $M_{1},\ldots,M_{n}$ be the
m-structure matrices of $A$ with respect to $B$. If $M_{i_{0}}$ is invertible
for some $1\leq i_{0}\leq n$ then $\mathrm{Ann}(A)=\{0\},$ and $A$ is an evolution
algebra if, and only if, each of the matrices $M_{i_{0}}^{-1}M_{1},...,M_{i_{0}}%
^{-1}M_{n}$ is diagonalisable (by similarity) for $j=1,...n$ and they pairwise commute.
\end{Corollary}

\begin{proof}
Since $\mathrm{Ann}(A)\subseteq\ker M_{i_{0}}$ by Lemma \ \ref{one}, we obtain
that if $M_{i_{0}}$ is invertible then $\mathrm{Ann}(A)=\{0\}$. Moreover, for
$\lambda_{0}=(0,...,0,\overset{(i_{0}-\text{th})}{1},0,...,0)$ we have%
\[
\mathrm{rank}(M(\lambda_{0}))=\mathrm{rank}(M_{i_{0}})=n
\]
and the result follows from Theorem \ref{rankn}.
\end{proof}
If $A$ is an algebra with $\mathrm{Ann}(A)\neq\{0\}$ (suppose that
$\dim\mathrm{Ann}(A)=r>0$) then we can fix a basis of\ $\mathrm{Ann}(A)$ which
can be extended to a basis of $A.$ Therefore we obtain a basis
$\,\widetilde{B}=\{e_{1},...,e_{r},e_{r+1},...,e_{n}\}$ of $A$ such that
$\{e_{r+1},...,e_{n}\}$ is a basis of $\mathrm{Ann}(A)$ and the m-structure
matrices $M_{1}(\widetilde{B}),\ldots,M_{n}(\widetilde{B})$ of $A$ with
respect to $\,\widetilde{B}$ satisfy $M_{k}(\widetilde{B}%
)=\widetilde{M}_{k}\oplus0_{n-r},$ for certain $r\times r$ matrices
$\widetilde{M}_{k}\in\mathcal{MS}_{r}.$

\begin{Theorem}
\label{dos}Let $A$ be a commutative complex algebra with $\mathrm{Ann}%
(A)\neq\{0\}.\,\ $Let $\,\widetilde{B}=\{e_{1},...,e_{r},e_{r+1},...,e_{n}\}$
be a basis of $A$ such that $\{e_{r+1},...,e_{n}\}$ is a basis of
$\mathrm{Ann}(A).$ Let $M_{1}(\widetilde{B}),...,M_{n}(\widetilde{B})$ be the
m-structure matrices of $A$ with respect to $\,\widetilde{B}$ with
$M_{k}(\widetilde{B})=\widetilde{M}_{k}\oplus0_{n-r},$ where $\widetilde{M}%
_{k}\in\mathcal{MS}_{r}.$ Then $A$ is an evolution algebra if, and only if,
there exists $\left\Vert \lambda_{0}\right\Vert =1$ such that the pencil
$\widetilde{M}(\lambda_{0})$ is invertible, each of the matrices
$\widetilde{M}(\lambda_{0})^{-1}\widetilde{M}_{1,}...,\widetilde{M}%
(\lambda_{0})^{-1}\widetilde{M}_{n,}$ is diagonalisable by similarity and they
pairwise commute.
\end{Theorem}

\begin{proof}
Assume $A$ is as stated. Then there exists $\left\Vert \lambda_{0}\right\Vert =1$ such that
the pencil $\widetilde{M}(\lambda_{0})$ is invertible if, and only if, the
maximum pencil rank of $M_{k}(\widetilde{B})$ is $r.$\ If this happens then
$\dim(\cap_{j=1}^{n}\ker M_{j}(\widetilde{B}))=n-r,$ as $\dim\mathrm{Ann}%
(A)=\dim(\cap_{j=1}^{n}\ker M_{j}(\widetilde{B}))$ by Lemma \ \ref{one}. Therefore if
$\widetilde{M}(\lambda_{0})$ is invertible then, by Corollary \ref{rankncoro},
we have that $\widetilde{M}_{1},...,\widetilde{M}_{n}$ are SDC if, and only
if, each of the matrices $\widetilde{M}(\lambda_{0})^{-1}\widetilde{M}%
_{1},...,$ $\widetilde{M}(\lambda_{0})^{-1}\widetilde{M}_{n}$ is
diagonalisable by similarity and they pairwise commute. Since the matrices
$\widetilde{M}_{1},...,\widetilde{M}_{n}$ are SDC (by $P_{r}\in\mathcal{GL}%
_{r}$) if, and only if, the matrices $M_{1}(\widetilde{B}),...,M_{n}%
(\widetilde{B})$ are SDC (by $P_{n}:=P_{r}\oplus I_{n-r}$), the result follows
from Theorem \ref{thm-main}.
\end{proof}

\begin{Remark}
\label{nota1}\emph{The above result shows that the condition that
}$A/\mathrm{Ann}(A)$\emph{ be an evolution algebra is a necessary condition
for }$A$ \emph{to be an evolution algebra. This is known because it was proved
in \cite{Ca-Si-Ve} that the quotient of an evolution algebra by an ideal is an
evolution algebra. However, Theorem \ref{dos} proves that this condition is not
sufficient (which is new). In fact, if }$\dim\mathrm{Ann}(A):=r<n,$\emph{ and
we consider a basis }$\,\widetilde{B},$ \emph{as in Theorem \ref{dos} above,
with m-structure matrices given by }$M_{k}(\widetilde{B})=\widetilde{M}%
_{k}\oplus0_{n-r}$\emph{ for }$k=1,...,n,$ \emph{then }$A$ \emph{is an
evolution algebra if, and only if, }$\widetilde{M}_{1,}...,\widetilde{M}_{n}%
$\emph{ are SDC. Suppose now that }$\widetilde{M}_{1,}...,\widetilde{M}_{r}%
$\emph{ are SDC but that }$\widetilde{M}_{1,}...,\widetilde{M}_{n}$\emph{ are
not SDC. It turns out that }$A/\mathrm{Ann}(A)$ \emph{is an evolution algebra
but }$A$\emph{ is not (because the m-structure matrices of }$A/\mathrm{Ann}%
(A)$ \emph{with respect to the basis }$\widetilde{B}_{A/\mathrm{Ann}%
(A)}=\{e_{1}+\mathrm{Ann}(A),...,e_{r}+\mathrm{Ann}(A)\}$\emph{ are precisely
}$\widetilde{M}_{1,}...,\widetilde{M}_{r}$\emph{). It is easy to come up with
particular examples of this situation (see Remark \ref{nota2} below).}
\end{Remark}

We conclude this section by providing a procedure, obtained from Theorems
\ref{thm-main}, \ \ref{rankn}, \ \ref{diago} and \ \ref{dos} above, to
determine in a finite number of steps whether or not a given commutative
algebra $A\,$ with fixed basis $B=\{e_{1},...,e_{n}\}$ is an evolution
algebra. Let $M_{1},...,M_{n}$ be the m-structure matrices of $A\,$ with
respect to $B.$

While one can try to check directly, see Example \ref{alg1} below, if the
matrices $M_{1},...,M_{n}$ are SDC this is generally not easy to do.
Alternatively, to determine if $A$ is an evolution algebra we can proceed as follows.

Check if any one of the matrices $M_{1},...,M_{n}$ is invertible.

(a) Suppose that $M_{i_{0}}$ is invertible, for some $1\leq i_{0}\leq n.$ If
$M_{i_{0}}^{-1}M_{1},...,M_{i_{0}}^{-1}M_{n}$ are all diagonalisable (by similarity) and they
pairwise commute then we can conclude that $A\,$\ is an evolution algebra, and
otherwise we conclude that $A$ is not an evolution algebra.

(b) If none of the matrices $M_{1},...,M_{n}$ is invertible then we determine
$Ann(A),$ that is, by means of (\ref{eq:M_k}), we describe those elements
$a\in A$ such that $ae_{i}=0\,$\ for every $i=1,...,n$.

(b.1) If $Ann(A)=\{0\}$ then we check if there exists some $\lambda
_{0}=(\lambda_{1},...,\lambda_{n})\in\mathbb{C}^{n}$ with $\left\Vert
\lambda_{0}\right\Vert =1$ such that $M(\lambda_{0}):=%
{\displaystyle\sum\limits_{i=1}^{n}}
\lambda_{i}M_{i}$ is invertible. If such a $\lambda_{0}$ does not exist then
we conclude that $A\,$ is not an evolution algebra. Otherwise we have that
$A$ is an evolution algebra if, and only if, the matrices $M(\lambda_{0}%
)^{-1}M_{1},...,M(\lambda_{0})^{-1}M_{n}$ are all diagonalisable (by similarity) and they
pairwise commute.

(b.2) If $Ann(A)\neq\{0\}$ then we construct a basis $\widetilde{B}%
=\{\widetilde{e}_{1},...,\widetilde{e}_{r},\widetilde{e}_{r+1}%
,...,\widetilde{e}_{n}\}$, such that $\{\widetilde{e}_{r+1},...,\widetilde{e}%
_{n}\}$ is a basis of $Ann(A)\neq\{0\}.$ We then have $M_{k}(\widetilde{B}%
)=\widetilde{M}_{k}\oplus0_{n-r}$ $\ $for $k=1,...,n$ and $r\times r$ matrices
$\widetilde{M}_{1},...,\widetilde{M}_{n}$. Next, we check if there exists
$\lambda_{0}=(\lambda_{1},...,\lambda_{n})\in\mathbb{C}^{n}$ with $\left\Vert
\lambda_{0}\right\Vert =1$ such that $\widetilde{M}(\lambda_{0}):=%
{\displaystyle\sum\limits_{i=1}^{n}}
\lambda_{i}\widetilde{M}_{i}$ is invertible as an $r\times r$ matrix. In
particular, this is the\textbf{ }case whenever $\widetilde{M}_{i_{0}}$ is
invertible for some $1\leq i_{0}\leq n$ (in which case we can choose
$\widetilde{M}(\lambda_{0})=\widetilde{M}_{i_{0}}$). If such a $\lambda_{0}$
does not exist then we conclude that $A$ is not an evolution algebra.
Otherwise, we have that $A$ is an evolution algebra if, and only if,\ the
matrices $\widetilde{M}(\lambda_{0})^{-1}\widetilde{M}_{1},...,$
$\widetilde{M}(\lambda_{0})^{-1}\widetilde{M}_{n}$ are all diagonalisable (by similarity) and
they pairwise commute.

\section{Some examples and applications}

We discuss some examples where our approach is useful to determine whether or
not certain classical genetic algebras are evolution algebras. Mostly these
algebras are defined in the literature as real algebras but, in our case, they
can be regarded as complex algebras (with the same basis, and hence with the
same m-structure matrices) as shown in Theorem \ref{complex1} and
Corollary \ref{complex2}.

We will consider the class of gametic algebras discussed by Etherington
\cite{etherington1940xxiii}.  Gametic algebras, widely used in genetics, are
simply baric algebras: they are endowed with a weight function. While further background is not necessary to decide if these algebras are evolution algebras or
not, we nevertheless refer the reader to
\cite{Liu} and \cite{Reed} for a review of these algebras.

\begin{Example}
\label{alg1}\emph{Let }$A$\emph{ be the algebra with basis }$B=\{e_{1}%
,e_{2}\}$\emph{ and }$e_{1}^{2}=e_{1}$\emph{, }$e_{1}e_{2}=e_{2}=e_{2}e_{1}%
,$\emph{ }$e_{2}^{2}=e_{1}.$\emph{ Define }$\xi:A\rightarrow\mathbb{K}%
~$\emph{by }$\xi(\alpha e_{1}+\beta e_{2})=\alpha+\beta.$\emph{ Obviously
}$\xi$\emph{ is linear and if }$a=\alpha e_{1}+\beta e_{2}$\emph{ and if
}$b=\gamma e_{1}+\delta e_{2}$\emph{ then}%
\[
ab=(\alpha\gamma+\beta\delta)e_{1}+(\alpha\delta+\beta\gamma)e_{2},
\]
\emph{so that }$\xi(ab)=(\alpha\gamma+\beta\delta)+(\alpha\delta+\beta
\gamma)=(\alpha+\beta)(\gamma+\delta)=\xi(a)\xi(b),$\emph{ and hence }$\xi
$\emph{ is a non-zero algebra homomorphism. Consequently }$A$\emph{ is a baric
algebra \cite{etherington1940xxiii}. }

\emph{The corresponding m-structure matrices with respect to }$B$\emph{ are
}$M_{1}=\left(
\begin{array}
[c]{cc}%
1 & 0\\
0 & 1
\end{array}
\right)  $\emph{ and }$M_{2}=\left(
\begin{array}
[c]{cc}%
0 & 1\\
1 & 0
\end{array}
\right)  .$ \emph{ Since for }$P=\left(
\begin{array}
[c]{rr}%
1 & 1\\
1 & -1
\end{array}
\right)  $ \emph{ we have that}$P^{T}M_{1}P=\left(
\begin{array}
[c]{cc}%
2 & 0\\
0 & 2
\end{array}
\right)  $ \emph{and} \ $P^{T}M_{2}P=\left(
\begin{array}
[c]{cc}%
2 & 0\\
0 & -2
\end{array}
\right)  ,$ \emph{by Theorem} \emph{\ref{thm-main}}, \emph{we obtain that }%
$A$\emph{ is an evolution algebra. In fact, }$\widetilde{B}=\{\widetilde{e}%
_{1},\widetilde{e}_{2}\}$,\emph{ with }$\widetilde{e}_{1}=e_{1}-e_{2}$\emph{
and }$\widetilde{e}_{2}=e_{1}+e_{2}$,\emph{ is a natural basis of }$A,$\emph{
as }$\widetilde{e}_{1}\widetilde{e}_{2}=0.$
\end{Example}

\begin{Remark}
\label{nota2}\emph{Let\ }$M_{1}$\emph{ and }$M_{2}$\emph{ be as above and
consider a matrix }$M_{3}$\emph{ that does not commute with }$M_{2}$\emph{,
say for instance }$M_{3}=\left(
\begin{array}
[c]{cc}%
1 & 0\\
0 & -1
\end{array}
\right)  .\ $\emph{Then we have that }$M_{1}^{-1}M_{2}$\emph{ and }$M_{1}%
^{-1}M_{3}$\emph{ do not commute so that, by the proof of Theorem \ref{dos}
(or alternatively using \cite[Section 3.3]{Bu-Me-Ve-SD}), the }$3\times3$
\emph{matrices }$M_{1}\oplus0_{1\times1},$\emph{ }$M_{2}\oplus0_{1\times1}%
$\emph{ and }$M_{3}\oplus0_{1\times1}$\emph{ are not SDC, while }$M_{1}$\emph{
and }$M_{2}$\emph{ are SDC. Therefore the algebra }$\widetilde{A}$\emph{ with
basis }$\widetilde{B}=\{e_{1},e_{2},e_{3}\}$\emph{ and product }$e_{1}%
^{2}=e_{1}+e_{3},$\emph{ }$\ e_{2}^{2}=e_{1}-e_{3},$\emph{ }$e_{3}^{2}%
=0,$\emph{ }$e_{1}e_{2}=e_{2}=e_{2}e_{1},$\emph{ }$e_{1}e_{3}=e_{3}e_{1}%
=e_{2}e_{3}=e_{3}e_{2}=0$\emph{ is an algebra such that }$\mathrm{Ann}%
(\widetilde{A})=\mathbb{K}e_{3}$\emph{. By Theorem \ref{dos} (see also Remark
\ref{nota1})}$\ $\emph{we have that }$\widetilde{A}$\emph{ is therefore not an
evolution algebra whereas }$\widetilde{A}/\mathrm{Ann}(\widetilde{A})$\emph{
is an evolution algebra isomorphic to the evolution algebra }$A$\emph{ in
Example \ref{alg1}.}
\end{Remark}

\begin{Example}
[Gametic algebra for simple Mendelian inheritance]\label{ejem1}

\emph{Let }$A_{0}$\emph{ denote a commutative }$2$\emph{-dimensional algebra
over\thinspace\ }$\mathbb{R}$\emph{, corresponding to the gametic algebra
describing simple Mendelian inheritance (see \cite{Reed}). In terms of the
basis \thinspace}$B=\{e_{1},e_{2}\}$\emph{ the multiplication table is }%
\[
e_{1}^{2}=e_{1},\quad e_{1}e_{2}=e_{2}e_{1}=\frac{1}{2}(e_{1}+e_{2}),\quad
e_{2}^{2}=e_{2}\,.
\]
\emph{The associated m-structure matrices }$M_{1},M_{2}$\emph{ can be read off
easily: }%
\[
M_{1}=\left(
\begin{array}
[c]{cc}%
1 & \frac{1}{2}\\
\frac{1}{2} & 0
\end{array}
\right)  ,\qquad M_{2}=\left(
\begin{array}
[c]{cc}%
0 & \frac{1}{2}\\
\frac{1}{2} & 1
\end{array}
\right)  \,.
\]
\emph{It is easy to check that }$A_{0}$\emph{ is a baric algebra, with weight
function defined by }$\xi(e_{1})=\xi(e_{2})=1$. \emph{Note that }$M_{1}%
^{-1}=\left(
\begin{array}
[c]{cc}%
0 & 2\\
2 & -4
\end{array}
\right)  $\emph{ while }%
\[
M_{1}^{-1}M_{2}=\emph{\ }\left(
\begin{array}
[c]{cc}%
0 & 2\\
2 & -4
\end{array}
\right)  \left(
\begin{array}
[c]{cc}%
0 & \frac{1}{2}\\
\frac{1}{2} & 1
\end{array}
\right)  =\allowbreak\left(
\begin{array}
[c]{cc}%
1 & 2\\
-2 & -3
\end{array}
\right)
\]
\emph{ is not diagonalisable by similarity, as }$\lambda=-1$\emph{ is the
unique eigenvalue and the associated eigenspace has dimension }$1.$\emph{
Therefore, by Corollary \ref{diago}, we obtain that }$A_{0}$\emph{ is not an
evolution algebra. (This last assertion can also be deduced from Theorem
\ref{thm-main}, with more tedious calculations, by directly checking that
}$M_{1}$ \emph{and }$M_{2}$\emph{ are not SDC).}
\end{Example}

We will now deform this algebra in order to construct an evolution algebra.

\begin{Example}
[Evolution algebra for deformed Mendelian inheritance]\label{ejem1de}

\emph{Consider a deformation of the algebra }$A_{0}$\emph{ of the previous
example. We denote these deformed algebras by }$A_{\varepsilon}$\emph{, which
depend on the free parameter }$\varepsilon\in\mathbb{R}$\emph{. In terms of
the basis }$B=\{e_{1},e_{2}\},$\emph{ the multiplication table for
}$A_{\varepsilon}$\emph{ is given by }%
\[
e_{1}^{2}=(1-\varepsilon)e_{1}+\varepsilon\,e_{2},\quad e_{1}e_{2}=e_{2}%
e_{1}=\frac{1}{2}(e_{1}+e_{2}),\quad e_{2}^{2}=e_{2}\,.
\]
\emph{The associated m-structure matrices }$M_{1},M_{2}$\emph{ are now: }%
\[
M_{1}=\left(
\begin{array}
[c]{cc}%
1-\varepsilon & \frac{1}{2}\\
\frac{1}{2} & 0
\end{array}
\right)  ,\qquad M_{2}=\left(
\begin{array}
[c]{cc}%
\varepsilon & \frac{1}{2}\\
\frac{1}{2} & 1
\end{array}
\right)  .
\]
\emph{For genetic applications we restrict }$0<\varepsilon\leq1$\emph{ so that
all coefficients in these matrices are non-negative. Moreover, }%
$A_{\varepsilon}$\emph{ is baric with weight function defined by }$\xi
(e_{1})=\xi(e_{2})=1$\emph{, for any }$\varepsilon$\emph{. In fact }$\xi
(e_{i}e_{j})=\xi(e_{i})\xi(e_{j})=1,$\emph{ for }$i,j=1,2$\emph{. Obviously,
the undeformed case corresponds to }$\varepsilon=0$.

\emph{Let us consider whether }$A_{\varepsilon}$\emph{ is an evolution algebra by
using Theorem \ref{rankn}. First of all, the maximal rank of the linear pencil
}$M(\mbox{\boldmath{$\lambda$}})=\lambda_{1}M_{1}+\lambda_{2}M_{2}$\emph{ is
}$r=2$\emph{ because }$M_{1}$\emph{ is nonsingular for all }$\varepsilon
$\emph{, so we can take }$\mbox{\boldmath{$\lambda$}}_{0}=(1,0)$\emph{. Thus
}$M(\mbox{\boldmath{$\lambda$}}_{0})=M_{1}$\emph{. To see that }%
$A_{\varepsilon}\,\ $\emph{is an evolution algebra we prove that }$M_{1}%
^{-1}M_{2}$\emph{ is diagonalisable by similarity. It is easy to check that}
\[
M_{1}^{-1}M_{2}=\left(
\begin{array}
[c]{cc}%
1 & 2\\
4\varepsilon-2 & 4\varepsilon-3
\end{array}
\right)
\]

\emph{and that if }%
\[
P=\left(
\begin{array}
[c]{cc}%
1 & 1\\
-1 & 2\varepsilon-1
\end{array}
\right)
\]
\emph{then }%

\begin{align*}
P^{-1}M_{1}^{-1}M_{2}P  &  =\\
&  =\left(
\begin{array}
[c]{cc}%
\frac{1}{2\varepsilon}\left(  2\varepsilon-1\right)  & -\frac{1}{2\varepsilon
}\\
\frac{1}{2\varepsilon} & \frac{1}{2\varepsilon}%
\end{array}
\right)  \left(
\begin{array}
[c]{cc}%
1 & 2\\
4\varepsilon-2 & 4\varepsilon-3
\end{array}
\right)  \left(
\begin{array}
[c]{cc}%
1 & 1\\
-1 & 2\varepsilon-1
\end{array}
\right) \\
&  =\left(
\begin{array}
[c]{cc}%
-1 & 0\\
0 & 4\varepsilon-1
\end{array}
\right)  .
\end{align*}

$\allowbreak$\emph{ Since}

\begin{align*}
P^{T}M(\mbox{\boldmath{$\lambda$}}_{0})P  &  =P^{T}M_{1}P=\left(
\begin{array}
[c]{cc}%
1 & -1\\
1 & 2\varepsilon-1
\end{array}
\right)  \left(
\begin{array}
[c]{cc}%
1-\varepsilon & \frac{1}{2}\\
\frac{1}{2} & 0
\end{array}
\right)  \left(
\begin{array}
[c]{cc}%
1 & 1\\
-1 & 2\varepsilon-1
\end{array}
\right) \\
&  =\left(
\begin{array}
[c]{cc}%
-\varepsilon & 0\\
0 & \varepsilon
\end{array}
\right)  ,
\end{align*}

\emph{and }$\det P=2\varepsilon$\emph{, we conclude by Theorem \ref{rankn}
that the algebra }$A_{\varepsilon}$\emph{ is an evolution algebra if, and only
if, }$\varepsilon\neq0$\emph{. For completeness we show the diagonalisation of
the original matrices: }%
\[
P^{T}M_{1}P=\left(
\begin{array}
[c]{cc}%
-\varepsilon & 0\\
0 & \varepsilon
\end{array}
\right)  ,\qquad P^{T}M_{2}P=\left(
\begin{array}
[c]{cc}%
\varepsilon & 0\\
0 & \varepsilon(4\varepsilon-1)
\end{array}
\right)  ,
\]
\emph{which shows by Theorem \ref{thm-main}, that }$A_{\varepsilon}$\emph{ is
an evolution algebra for every }$\varepsilon>0,$ \emph{having }$B=\{e_{1}%
-e_{2},$ $e_{1}+(2\varepsilon-1)e_{2}\}$\emph{ as a natural basis.}
\end{Example}

\begin{Example}
\emph{The annihilator of every algebra }$A_{\varepsilon}$\emph{ in the above
example is zero as one of \ its m-structure matrices is invertible. To get a
similar example with algebras having non-zero annihilator, consider for
instance the algebra }$A_{\varepsilon}$\emph{ with natural basis }%
$\widehat{B}=\{e_{1},e_{2},e_{3}\}$\emph{ and product given by }%

\begin{align*}
e_{1}^{2}  &  =(1-\varepsilon)e_{1}+\varepsilon\,e_{2}-\varepsilon
\,e_{3},\text{ \ }e_{2}^{2}=e_{2}-e_{3};\text{ }e_{3}^{2}=0,\\
e_{1}e_{2}  &  =e_{2}e_{1}=\frac{1}{2}(e_{1}+e_{2}-e_{3}),\text{ \quad}%
e_{1}e_{3}=e_{3}e_{1}=e_{2}e_{3}=e_{3}e_{2}=0.\,
\end{align*}

\emph{Here the m-structure matrices are }$M_{k}(\widehat{B})=M_{k}\oplus
0$\emph{ (for }$i=1,2,3)$,\emph{ where }$0$\emph{ denotes the }$1\times1$\emph{
zero matrix, }$M_{1}$\emph{ and }$M_{2}$\emph{ are given in the above example
and }$M_{3}=-M_{2}.$\emph{ Hence if }%
\[
P=\left(
\begin{array}
[c]{ccc}%
1 & 1 & 0\\
-1 & 2\varepsilon-1 & 0\\
0 & 0 & 1
\end{array}
\right)
\]
\emph{we obtain, from the calculations in the above example, that }$P^{T}%
M_{k}(\widehat{B})P$\emph{ is diagonal for every }$k=1,2,3$\emph{ and hence
}$A_{\varepsilon}$\emph{ is an evolution algebra. Nevertheless, for
}$\varepsilon=0$\emph{ we do not obtain an evolution algebra. Indeed, if we
denote this algebra by }$A$\emph{ then the quotient algebra }$A/\mathrm{Ann}%
(A)$\emph{ is exactly the algebra }$A_{0}$\emph{ in Example \ref{ejem1} which
is not an evolution algebra and, consequently, }$A$\emph{ is not an evolution
algebra (see Remark \ref{nota1}).}
\end{Example}

\begin{Example}
[Gametic algebra for auto-tetraploid inheritance]\label{ejem2}

\emph{Let }$T_{0}$\emph{ denote a }$3$\emph{-dimensional commutative algebra
over }$\mathbb{R}$\emph{, considered the simplest case of special train
algebras in polyploidy \cite[Chapter 15]{etherington1940xxiii} (see also
\cite{Liu} and \cite{Reed}). In terms of the basis }$\{e_{1},e_{2},e_{3}%
\}$\emph{ the multiplication table is given by }

\begin{align*}
e_{1}^{2}  &  =e_{1},\quad e_{2}^{2}=e_{1}e_{3}=\frac{1}{6}(e_{1}+4e_{2}%
+e_{3}),\quad\\
e_{3}^{2}  &  =e_{3},\quad e_{2}e_{3}=\frac{1}{2}(e_{2}+e_{3}),\quad
e_{1}e_{2}=\frac{1}{2}(e_{1}+e_{2})\,.
\end{align*}

\emph{The corresponding m-structure matrices }$M_{1},M_{2},M_{3}$\emph{ are }%
\[
M_{1}=\left(
\begin{array}
[c]{ccc}%
1 & \frac{1}{2} & \frac{1}{6}\\
\frac{1}{2}\smallskip & \frac{1}{6} & 0\\
\frac{1}{6} & 0 & 0
\end{array}
\right)  ,\quad M_{2}=\left(
\begin{array}
[c]{ccc}%
0 & \frac{1}{2} & \frac{2}{3}\\
\frac{1}{2}\smallskip & \frac{2}{3} & \frac{1}{2}\\
\frac{2}{3} & \frac{1}{2} & 0
\end{array}
\right)  ,\quad M_{3}=\left(
\begin{array}
[c]{ccc}%
0 & 0 & \frac{1}{6}\\
0\smallskip & \frac{1}{6} & \frac{1}{2}\\
\frac{1}{6} & \frac{1}{2} & 1
\end{array}
\right)  .
\]
\emph{The algebra }$T_{0}$\emph{ is baric, with weight function defined by
}$\xi(e_{j})=1,\,\,\,j=1,2,3$.

\emph{To see that this algebra is not an evolution algebra note that }%
\[
M_{1}^{-1}=\left(
\begin{array}
[c]{rrr}%
0 & 0 & 6\\
0 & 6 & -18\\
6 & -18 & 18
\end{array}
\right)  ,
\]
\emph{and that}%
\[
M_{1}^{-1}M_{2}=\left(
\begin{array}
[c]{ccc}%
0 & 0 & 6\\
0 & 6 & -18\\
6 & -18 & 18
\end{array}
\right)  \left(
\begin{array}
[c]{ccc}%
0 & \frac{1}{2} & \frac{2}{3}\\
\frac{1}{2} & \frac{2}{3} & \frac{1}{2}\\
\frac{2}{3} & \frac{1}{2} & 0
\end{array}
\right)  =\left(
\begin{array}
[c]{rrr}%
4 & 3 & 0\\
-9 & -5 & 3\\
3 & 0 & -5
\end{array}
\right)
\]
\emph{ is not diagonalisable by similarity because it has a single eigenvalue
(}$\lambda=-2$\emph{) and the dimension of the associated eigenspace is }%
$1$\emph{ (indeed, }$(1,-2,1)^{T}$ \emph{generates it)}$.$\emph{ Consequently,
}$A$ \emph{is not an evolution algebra} \emph{by Corollary \ref{diago}. }

\emph{On the other hand, }%
\[
M_{1}^{-1}M_{3}=\left(
\begin{array}
[c]{ccc}%
0 & 0 & 6\\
0 & 6 & -18\\
6 & -18 & 18
\end{array}
\right)  \left(
\begin{array}
[c]{ccc}%
0 & 0 & \frac{1}{6}\\
0 & \frac{1}{6} & \frac{1}{2}\\
\frac{1}{6} & \frac{1}{2} & 1
\end{array}
\right)  =\left(
\begin{array}
[c]{rrr}%
1 & 3 & 6\\
-3 & -8 & -15\\
3 & 6 & 10
\end{array}
\right)
\]
\emph{so that }%
\[
M_{1}^{-1}M_{2}M_{1}^{-1}M_{3}=M_{1}^{-1}M_{3}M_{1}^{-1}M_{2}=\left(
\begin{array}
[c]{rrr}%
-5 & -12 & -21\\
15 & 31 & 51\\
-12 & -21 & -32
\end{array}
\right)  .
\]
\emph{This proves that, in Theorem \ref{rankn}, the condition that the
matrices }$M(\lambda_{0})^{-1}M_{1},...,$ $M(\lambda_{0})^{-1}M_{n}$
\emph{\ pairwise commute is not sufficient to ensure that the given algebra is
an evolution algebra (see also Proposition \ref{Horn}).}
\end{Example}

\begin{Example}
[Evolution algebra for deformed auto-tetraploid inheritance]\label{ejem2de}

\emph{Consider now a deformation of the algebra }$T_{0}$\emph{ of the previous
example. We denote this deformed algebra by }$T_{\varepsilon}$\emph{, which
depends on the free parameter }$\varepsilon\in\mathbb{R}$\emph{. In terms of
the basis }$\{e_{1},e_{2},e_{3}\}$\emph{ the multiplication table for
}$T_{\varepsilon}$\emph{ is: }%
\[
e_{1}^{2}=e_{1}+2\,\varepsilon(e_{1}+4e_{2}),\quad e_{2}^{2}=\frac{1}{6}%
(e_{1}+4e_{2}+e_{3})-\varepsilon(3e_{2}-13e_{3}),\quad e_{3}^{2}%
=e_{3}+10\,\varepsilon e_{3}\,,
\]%
\[
e_{1}e_{3}=\frac{1}{6}(e_{1}+4e_{2}+e_{3})+10\,\varepsilon e_{3},\text{ }%
e_{2}e_{3}=\frac{1}{2}(e_{2}+e_{3})+10\,\varepsilon e_{3},\text{ }e_{1}%
e_{2}=\frac{1}{2}(e_{1}+e_{2})+10\,\varepsilon e_{3}\,.
\]
\emph{The corresponding m-strucuture matrices }$M_{1},M_{2},M_{3}$\emph{ are }%
\[
M_{1}=\left(
\begin{array}
[c]{ccc}%
1+2\,\varepsilon & \frac{1}{2} & \frac{1}{6}\\
\frac{1}{2} & \frac{1}{6} & 0\\
\frac{1}{6} & 0 & 0
\end{array}
\right)  ,\quad M_{2}=\left(
\begin{array}
[c]{ccc}%
8\,\varepsilon & \frac{1}{2} & \frac{2}{3}\\
\frac{1}{2} & \frac{2}{3}-3\,\varepsilon & \frac{1}{2}\\
\frac{2}{3} & \frac{1}{2} & 0
\end{array}
\right)  ,
\]%
\[
M_{3}=\left(
\begin{array}
[c]{ccc}%
0 & 10\,\varepsilon & \frac{1}{6}+10\,\varepsilon\\
10\,\varepsilon & \frac{1}{6}+13\,\varepsilon & \frac{1}{2}+10\,\varepsilon\\
\frac{1}{6}+10\,\varepsilon & \frac{1}{2}+10\,\varepsilon & 1+10\,\varepsilon
\end{array}
\right)  .
\]
\emph{For genetic applications, we restrict }$0<\varepsilon\leq2/9$\emph{, so
all coefficients in the above matrices are non-negative. The algebra
}$T_{\varepsilon}$\emph{ is baric, with weight function defined by }$\xi
(e_{j})=1+10\,\varepsilon,\,\,\,j=1,2,3$.

\emph{Let us consider whether }$T_{\varepsilon}$\emph{ is an evolution
algebra. }$\ $\emph{First of all, the maximal rank of the linear pencil
}$M(\mbox{\boldmath{$\lambda$}})=\lambda_{1}M_{1}+\lambda_{2}M_{2}+\lambda
_{3}M_{3}$\emph{ is }$r=3$\emph{ because }$M_{1}$\emph{ is nonsingular for all
}$\varepsilon$\emph{, so we can take }$\mbox{\boldmath{$\lambda$}}_{0}%
=(1,0,0)$\emph{. Thus }$M(\mbox{\boldmath{$\lambda$}}_{0})=M_{1}$\emph{. By
Theorem \ref{rankn}, a necessary condition is that the matrices }$M_{1}%
^{-1}M_{2}$\emph{ and }$M_{1}^{-1}M_{3}$\emph{ are simultaneously
diagonalisable by similarity: in particular, they must commute. Let us write
these matrices explicitly: }%
\[
M_{1}^{-1}M_{2}=\left(
\begin{array}
[c]{ccc}%
4 & 3 & 0\\
-9 & -5-18\,\varepsilon & 3\\
3 & 18\,\varepsilon & -5
\end{array}
\right)  \,,
\]%
\[
M_{1}^{-1}M_{3}=\left(
\begin{array}
[c]{ccc}%
1+60\,\varepsilon & 3+60\,\varepsilon & 6+60\,\varepsilon\\
-3(1+40\,\varepsilon) & -2(4+51\,\varepsilon) & -15(1+8\,\varepsilon)\\
3(1-4\,\varepsilon-240\,\varepsilon^{2}) & 6(1-5\,\varepsilon-120\,\varepsilon
^{2}) & 2(5-6\,\varepsilon-360\,\varepsilon^{2})
\end{array}
\right)  \,.
\]
\emph{It is straightforward to show that these matrices commute for all
}$\varepsilon$\emph{ (even for }$\varepsilon=0$\emph{). Regarding the Jordan
decomposition for }$M_{1}^{-1}M_{2}$\emph{ and }$M_{1}^{-1}M_{3}$\emph{ we
find that if }$\varepsilon>0$\emph{ then these matrices are simultaneously
diagonalisable: in fact, there is a nonsingular matrix }$P$\emph{ such that
}$P^{-1}M_{1}^{-1}M_{2}P$\emph{ is diagonal: }%
\[
P^{-1}M_{1}^{-1}M_{2}P=\left(
\begin{array}
[c]{ccc}%
-2 & 0 & 0\\
0 & -2-9\varepsilon-3S_{\varepsilon} & 0\\
0 & 0 & 2-9\varepsilon+3S_{\varepsilon}%
\end{array}
\right)  ,\qquad S_{\varepsilon}=\sqrt{3\varepsilon(3\varepsilon+4)}\,.
\]
\emph{Explicitly, in terms of the radical }$S_{\varepsilon}$\emph{, }%
\[
P=\left(
\begin{array}
[c]{ccc}%
1 & 1 & 1\\
-2 & -2-3\varepsilon-S_{\varepsilon} & -2-3\varepsilon+S_{\varepsilon}\\
1-12\varepsilon & 1+3\varepsilon+S_{\varepsilon} & 1+3\varepsilon
-S_{\varepsilon}%
\end{array}
\right)  \,.
\]
\emph{We find }$\det P=-24\varepsilon S_{\varepsilon}$\emph{ which shows there
is a problem at }$\varepsilon=0$\emph{. It is easy to show that at
}$\varepsilon=0$\emph{ the Jordan form of }$M_{1}^{-1}M_{2}$\emph{ is not
diagonal. For }$\varepsilon>0$\emph{ the Jordan form of }$M_{1}^{-1}M_{2}%
$\emph{ is diagonal and so is the Jordan form of }$M_{1}^{-1}M_{3}$\emph{: }%
\[
P^{-1}M_{1}^{-1}M_{3}P=\left(
\begin{array}
[c]{ccc}%
1-72\varepsilon-720\varepsilon^{2} & 0 & 0\\
0 & 1+9\varepsilon+3S_{\varepsilon} & 0\\
0 & 0 & 1+9\varepsilon-3S_{\varepsilon}%
\end{array}
\right)  .
\]

\emph{For completeness we show the diagonalisation of the original matrices: }%
\[
P^{T}M_{1}P=\varepsilon\left(
\begin{array}
[c]{ccc}%
-2 & 0 & 0\\
0 & 4+3\varepsilon+S_{\varepsilon} & 0\\
0 & 0 & 4+3\varepsilon-S_{\varepsilon}%
\end{array}
\right)  ,
\]%
\[
{\small P}^{T}{\small M}_{2}{\small P=-2\varepsilon\,}\left(
\begin{array}
[c]{ccc}%
-2 & 0 & 0\\
0 & 4+39\varepsilon+27\varepsilon^{2}+(9\varepsilon+7)S_{\varepsilon} & 0\\
0 & 0 & 4+39\varepsilon+27\varepsilon^{2}-(9\varepsilon+7)S_{\varepsilon}%
\end{array}
\right)  ,
\]%
\[
{\small P}^{T}{\small M}_{3}{\small P=\varepsilon\,}\left(
\begin{array}
[c]{ccc}%
\alpha & 0 & 0\\
0 & {\small 4+75\varepsilon+54\varepsilon}^{2}{\small +(18\varepsilon
+13)S}_{\varepsilon} & 0\\
0 & 0 & {\small 4+75\varepsilon+54\varepsilon}^{2}{\small -(18\varepsilon
+13)S}_{\varepsilon}%
\end{array}
\right)  {\small ,}%
\]
\emph{where }$\alpha=-2+144\varepsilon+1440\varepsilon^{{\small 2}}.$
\end{Example}

\section{Conclusions and Discussion}

In this paper we determine completely whether a given algebra $A$ is an
evolution algebra, by translating the question to a recently solved problem,
namely, the problem of simultaneous diagonalisation via congruence of the
m-structure matrices of $A$. This is relevant because evolution algebras have
strong connections with areas such as group theory, Markov processes, theory of
knots, and graph theory, amongst others. In fact, every evolution algebra can be
canonically regarded as a weighted digraph when a natural basis is fixed, and
because of this evolution algebras may  introduce useful algebraic techniques
into the study of some digraphs.\textbf{ }

We also consider applications of our results to classical genetic algebras.
Strikingly, the classical cases of Mendelian and auto-tetraploid inheritance
are not evolution algebras, while slight deformations of them produce
evolution algebras. This is interesting because evolution algebras are
supposed to describe asexual reproduction, unlike these classical cases. In
future work we will study more closely the relation between baric algebras and
evolution algebras, in order to better understand this phenomenon.

\vspace{6pt}

\authorcontributions{All authors contributed equally to this manuscript}

\funding{This work was partially supported by Project MTM216-76327-C3-2-P .\\
This work was also supported by the award of the Distinguished Visitor Grant
of the School of Mathematics and Statistics, University College Dublin to the
third author}


\conflictsofinterest{The authors declare no conflict of interest.\\
\\
\textbf{Mathematics Subject Classification }[2010]: Primary 17D92 and 15A60.}

\reftitle{References}



\begin{thebibliography}{99}                                                                                               %
\bibitem {Tian}Tian J.P., \textit{Evolution algebras and their applications}.
Lecture Notes in Mathematics vol. \textbf{1921}. Springer-Verlag (2008).

\bibitem {Be-Me-Ve}Becerra J., Beltr\'{a}n M. and Velasco M.V., Pulse
Processes in Networks and Evolution Algebras,\textit{ Mathematics }\textbf{8
}(2020), 387-407.

\bibitem {Ca-Si-Ve}Cabrera Y., Siles M., and \ Velasco M.V., Evolution
algebras of arbitrary dimension and their decompositions, \textit{Linear
Algebra and its Applications }\textbf{495} (2016), 122-162.

\bibitem {Ce-Ve}Celorrio M. E. and Velasco M.V, Classifying evolution algebras
of dimension two and three,\textit{ Mathematics }\textbf{7} (2019), 1236-1261.

\bibitem {Ca-Go}Camacho L.M., G\'{o}mez, J.R., Omirov B.A. and Turdibaev R.M.,
Some properties of evolution algebras,\textit{ Bulletin of the Korean
Mathematical Society }\textbf{50} (2013), 1481--1494.

\bibitem {3}Camacho L.M., G\'{o}mez, J.R., Omirov B.A. and Turdibaev R.M., The
derivations of some evolution algebras,.\textit{ Linear Multilinear Algebra}
\textbf{6} (2013), 309--322.

\bibitem {5}Casas J.M., Ladra M. and Rozikov U.A., A chain of evolution
algebras, \textit{Linear Algebra and its Applications} \textbf{435} (2011), 852--870.

\bibitem {Elduque}Elduque A. and Labra A., On nilpotent evolution algebras,
\textit{Linear Algebra and its Applications} \textbf{505} (2016), 11-31.

\bibitem {He}Hegazi A.S. and Abdelwahab H., Nilpotent evolution algebras over
arbitrary fields, \textit{Linear Algebra and Applications} \textbf{486}
(2015), 345-360.

\bibitem {13}Labra A., Ladra M. and Rozikov U.A., An evolution algebra in
population genetics, .\textit{Linear Algebra and Applications} \textbf{45}
(2014), 348--362.

\bibitem {15}Ladra M., Omirov B.A. and Rozikov U.A., Dibaric and evolution
algebras in Biology, \textit{Lobachevskii Journal of Mathematics} \textbf{35}
(2014), 198--210.

\bibitem {Mellon-Ve}Mellon P. and Velasco M.V., Analytic aspects of evolution
algebras, \textit{Banach Journal of Mathematical Analysis} \textbf{13} (2019), 113-132.

\bibitem {Ro-Ve}Omirov B., Rozikov U. and Velasco M.V., A class of nilpotent
evolution algebras, \textit{Communications in Algebra } \textbf{47} (2019),
1556- 1567.

\bibitem {19}Rozikov U.A. and Tian J.P., Evolution algebras generated by Gibbs
measures,  \textit{Lobachevskii Journal of Mathematics}\textbf{ 32} (2011), 270--277.

\bibitem {20}Rozikov U.A. and Murodov Sh.N., Dynamics of two-dimensional
evolution algebras, \textit{Lobachevskii Journal of Mathematics} \textbf{3}
(2013), 344--358.

\bibitem {21}Rozikov U.A. and Velasco M.V., Discrete-time dynamical system and
an evolution algebra of mosquito population, \textit{Journal of Mathematical
Biology} \textbf{78} (2019), 1225-1244.

\bibitem {Ve}Velasco M.V., The Jacobson radical of an evolution algebra,
\textit{Journal of Spectral Theory}\textbf{ 9} (2019), 601-634.

\bibitem {Hi-U}Hiriart-Urruty J.B., Potpourri of conjectures and open
questions in nonlinear analysis and optimization, \textit{SIAM Review
}\textbf{49} (2007), 255--273.

\bibitem {Hi-Ur-Tor}Hiriart-Urruty J.B. and Torki M., Permanently going back
and forth between the \textquotedblleft quadratic world\textquotedblright\ and
the \textquotedblleft convexity world\textquotedblright\ in Optimization,
\textit{Applied Mathematics and} \textit{Optimization} \textbf{45} (2002), 169--184.

\bibitem {Ji-Li}Jiang R. and Li \ D., Simultaneous diagonalization of matrices
and its applications in quadratically constrained quadratic programming,
\textit{SIAM \ Journal on Optimization} \textbf{26 }(2016), 1649--1669.

\bibitem {belouchrani1997blind}Belouchrani A., Abed-Meraim K., Cardoso J.F.
and Moulines E., A blind source separation technique using second-order
statistics, IEEE Transactions on signal processing \textbf{45} (1997), 434--444.

\bibitem {Vol}Vollgraf R. and Obermayer K., Quadratic optimization for
simultaneous matrix diagonalization, \textit{IEEE Transactions on Signal
Processing} \textbf{54} (2006), 3270--3278.

\bibitem {yeredor2000blind}Yeredor A., Blind source separation via the second
characteristic function, \textit{Signal Processing } \textbf{80} (2000), 897--902.

\bibitem {yeredor2002non}Yeredor A., Non-orthogonal joint diagonalization in
the least-squares sense with application in blind source separation, \textit{IEEE Transactions on signal processing} \textbf{50} (2002), 1545--1553.

\bibitem {Bu-Me-Ve-SD} Bustamante M. D., Mellon P. and Velasco M.V., Solving the Problem of Simultaneous Diagonalization of Complex Symmetric Matrices via Congruence, SIAM J. Matrix Anal.  Appl., \textbf{41} (2020), 1616--1629.

\bibitem {Pierce82}Pierce R.S., \textit{Associative Algebras,} Springer, New
York, NY (1982).

\bibitem {Horn}Horn R.A. and Johnson C.R., \textit{Matrix Analysis}, Second
Edition. Cambridge University Press (2013).

\bibitem {etherington1940xxiii}Etherington I. M. H., XXIII.---Genetic
Algebras, \textit{Proceedings of the Royal Society of Edinburgh }\textbf{59}
(1940), 242--258.

\bibitem {Liu}Liubich Y.I.,.\textit{Mathematical Structures in Population
Genetics, }Springer-Verlag (1992).

\bibitem {Reed}Reed M., Algebraic structure of genetic inheritance,
\textit{Bulletin of the American Mathematical Society} \ \textbf{34} (1997), 107--130.

\end{thebibliography}





\end{document}